\documentclass{article}
\usepackage[latin1]{inputenc}
\usepackage[english]{babel}
\usepackage{amsfonts}
\usepackage{amssymb}
\usepackage{amsmath}
\usepackage{latexsym}
\newtheorem{theoreme}{Theorem}[section]
\newtheorem{lemme}[theoreme]{Lemma}
\newtheorem{corollaire}[theoreme]{Corollary}
\newenvironment{preuve}{\noindent\emph{Proof} : }{\begin{flushright}$\Box$\end{flushright}}
\newtheorem{remarque}[theoreme]{Remark}
\newtheorem{proposition}[theoreme]{Proposition}
\newtheorem{defi}[theoreme]{Définition}
\newcommand{\F}{\mathbb{F}}
\newcommand{\pgcd}{\mathrm{gcd}}
\begin{document}
\renewcommand{\labelitemi}{$\cdot$}
\title{Functions which are PN on infinitely many extensions of $\F_p$, $p$ odd.}
\author{Elodie Leducq}
\date{}
\maketitle

\hrule \vspace{0.5cm}
\noindent\textbf{\textit{Abstract}} : Jedlicka, Hernando and McGuire have proved that Gold and Kasami functions are the only power mappings which are APN on infinitely many extensions of $\F_2$. For $p$ an odd prime, we prove that the only power mappings $x\mapsto x^m$ such that $m\equiv1\mod p$ which are PN on infinitely many extensions of of $\F_p$ are those such that $m=1+p^l$, l positive integer. As Jedlicka, Hernando and McGuire, we prove that $\frac{(x+1)^m-x^m-(y+1)^m+y^m}{x-y}$ has an absolutely irreducible factor by using Bézout's Theorem.
\\
\hrule
\section{Introduction}

In \cite{jedlicka_excep} and \cite{mcguire_excep}, the authors are interested in integers $m$ such that the function $x\mapsto x^m$ is almost perfectly nonlinear (APN) on infinitely many extensions of $\F_2$. Using similar methods, we study here the case of perfectly nonlinear (PN) functions over finite fields of odd characteristic.
\\\\\indent Let $p$ be a prime number, $n$ a positive integer, $q=p^n$ and $\F_{q}$ a finite field with $q$ elements. We recall the following definition :
\begin{defi}We say that the function $\phi$ is APN over $\F_{q}$ if :
$$\forall a,b\in\F_{q},a\neq0,|\{x\in\F_{q},\phi(x+a)-\phi(x)=b\}|\leq2$$
and if, furthermore, there exists a pair $(a,b)$ such that we have equality.\end{defi}

Jedlicka, Hernando and McGuire prove that in the case of characteristic 2, the only integers $m$ such that $x\mapsto x^m$ is APN on infinitely many extensions of $\F_2$ are $m=2^k+1$ (Gold) and $m=4^k-2^k+1$ (Kasami). They use the fact that a function $x\mapsto x^m$ is APN over $\F_{2^n}$ if and only if the rational points in $\F_{2^n}$ of $(x+1)^m+x^m+(y+1)^m+y^m=0$ are points such that $x=y$ or $x=y+1$. This can happen only if $\frac{(x+1)^m+x^m+(y+1)^m+y^m}{(x+y)(x+y+1)}$ has no absolutely irreducible factor over $\F_2$.
\\\\\indent In characteristic 2, if $\phi(x+a)+\phi(x)=b$ then $\phi(x+a+a)+\phi(x+a)=b$. So, there is no PN function (see definition below). Furthermore, to prove that $\phi$ is APN, it is sufficient to prove that $$\forall a,b\in\F_{q},a\neq0,|\{x\in\F_{q},\phi(x+a)-\phi(x)=b\}|\leq2.$$
\indent On the contrary, in odd characteristic, we do not know the relationship between the two solutions of $\phi(x+a)-\phi(x)=b$ (if there exist two). So, it seems to be difficult to adapt this method to APN functions in odd characteristic. Nevertheless, we can try on PN functions :
\begin{defi}If $q$ is odd, a function $\phi$ is PN over $\F_{q}$ if for all $b\in\F_{q}$ and all $a\in\F_{q}^*$
$$|\{x\in\F_{q},\phi(x+a)-\phi(x)=b\}|=1.$$
Equivalently, a function $\phi$ is PN over $\F_{q}$ if for all $a\in\F_{q}^*$, the only rational points in $\F_{q}$ of
$$\phi(x+a)-\phi(x)-\phi(y+a)+\phi(y)=0$$ are points such that $x=y$.\end{defi}
From now, $\phi:x\mapsto x^m$, $m\geq3$. We only have to consider the case where $a=1$ in the definition above (see \cite{DOB_apn_odd_char}).
\begin{remarque}If $m$ is odd then, 0 and 1 are solutions of $(x+1)^m-x^m=1$. So $x \mapsto x^m$ is not PN over $\F_{p^n}$ for any $n$. \end{remarque}

The only known PN power mappings are the following :

\begin{proposition}\label{PNP}Let $f : x\mapsto x^m$ a power mapping. Then $f$ is PN on $\F_{p^n}$ for 
\begin{enumerate}\item $m=2$,
\item $m=p^l+1$ where $l$ is an integer such that $\frac{n}{\gcd(n,l)}$ is odd \cite{MR1432296,MR0226486},
\item $m=\frac{3^l+1}{2}$ where $p=3$ and $l$ is an odd integer such that $\gcd(l,n)=1$ \cite{MR1432296}.\end{enumerate}\end{proposition}

We set $f(x,y)=(x+1)^m-x^m-(y+1)^m+y^m$. Since $(x-y)$ divides $f(x,y)$, we define $h(x,y)=\frac{f(x,y)}{(x-y)}$.
\\\\\indent We can assume that $m\not\equiv0\mod p$. Indeed, if $x\mapsto x^m$ is PN over $\F_q$ and $m\equiv0\mod p$ then $x\mapsto x^{\frac{m}{p}}$ is also PN over $\F_q$.
\begin{proposition}\label{PN}If $h$ has an absolutely irreducible factor over $\F_p$ then $x\mapsto x^m$ is not PN on $\F_{p^n}$ for $n$ sufficiently large.\end{proposition}
\begin{preuve} Assume that $h$ has an absolutely irreducible factor over $\F_p$, denoted by $Q$. If $Q(x,y)=c(x-y)$ with $c\in\F_p^*$, then $f(x,y)=(y-x)^2\widetilde{Q}(x,y)$, $\widetilde{Q}\in\F_p[x,y]$.
Hence,
$$-m(y+1)^{m-1}+my^{m-1}=\frac{\partial f}{\partial y}(x,y)=2(y-x)\widetilde{Q}(x,y)+(y-x)^2\frac{\partial\widetilde{Q}}{\partial y}(x,y).$$
So, we get that for all $x\in\F_{p^n}$, $-m(x+1)^{m-1}+mx^{m-1}=0$ which is impossible since $m\not\equiv0\mod p$.
Let $s$ be the degree of $Q$. Since $Q\neq  c(x-y)$, $Q(x,x)$ is not the null polynomial. So there are at most $s$ rational points of $Q$ such that $x=y$.
\\\indent On the other hand, if we denote by P the number of affine rational points of $Q$ on $\F_{p^n}$, we have (see \cite[p. 331]{lidl_finite_fields}):
$$|P-p^n|\leq (s-1)(s-2)\sqrt{p^n}+s^2.$$ 
Hence, for $n$ sufficiently large, $Q$ has a rational point in $\F_{p^n}$ such that $x\neq y$ and $x\mapsto x^m$ is not PN over $\F_{p^n}$.\end{preuve}

From now, we are interested in the case where $m\equiv1\mod p$. We denote by $l$ the greatest integer such that $p^l$ divides $m-1$ and we set $$d:=\pgcd(m-1,p^l-1)=\pgcd\left(\frac{m-1}{p^l},p^l-1\right).$$
By Proposition \ref{PNP}, the only known functions $x\mapsto x^m$ such that $m\equiv 1\mod p$ which are PN on infinitely many extensions of $\F_p$ are those such that $m=1+p^l$. We want to prove that there does not exist any other.

\begin{theoreme}\label{facteurs} Let $m$ be an integer such that $m\geq3$, $m\equiv1\mod p$ and \\$m\neq 1+p^l$. Assume that $\frac{m-1}{p^l}\neq p^l-1$. Then $h$ has an absolutely irreducible factor over $\F_p$.\end{theoreme}

\begin{corollaire}The only $m\equiv1\mod p$ such that $x\mapsto x^m$ is PN on infinitely many extensions of $\F_p$ are $m=1+p^l$.\end{corollaire}

\begin{preuve}By Theorem \ref{facteurs} and Proposition \ref{PN}, we only have to treat the case where $d=\frac{m-1}{p^l}=p^l-1$. Then $m=p^l(p^l-1)+1$ which is odd; so $x\mapsto x^m$ is not PN on all extensions of $\F_p$\end{preuve}

\begin{remarque}After reading the manuscript of this paper, McGuire informed me that the result in this paper has already been proved by Hernando and himself and also by R. Coulter. But as far as I know, it has never been published.\end{remarque}

Now, we only have to prove Theorem \ref{facteurs}. The method of Jedlicka, Hernando and McGuire is, using Bézout's Theorem, to prove that $h$ has an absolutely irreducible factor over $\F_p$ because it has not enough singular points. In Part 2, we study singular points of $h$ and their multiplicity. In Part 3, we bound the intersection number (see \cite{fulton_alg_curves} for definition) $I_t(u,v)$ where $t$ is a singular point of $h$ and $u$, $v$ are such that $h=uv$. In part 4, we prove Theorem \ref{facteurs}. Finally, we consider briefly the case where $d\not\equiv1\mod p$.

\section{Singularities of $h$}\label{P2}
\begin{proposition}The singular points of $h$ are described in Table \ref{tab1}.\end{proposition}
 The proof of this theorem follows from Lemmas \ref{infini} to \ref{IIc} and their corollaries.

\begin{table}[h]
\hrule
\caption{\label{tab1}Singularities of $h$ for $m\equiv1\mod p$}
\begin{minipage}[b]{\textwidth}
\renewcommand{\footnoterule}{}
\begin{tabular}{|c|p{2.5cm}|c|c|c|}
\hline
Type & Description & $m_t(h)$ & $I_t$ bound & max number of points\\
\hline
Ia & Affine $x_0=y_0$ $x_0$, $y_0\in\F_{p^l}^*$ & $p^l$ & $\frac{p^{2l}-1}{4}$ & $d-1$\\
\hline
Ib & Affine $x_0=y_0$, $x_0$, $y_0\not\in\F_{p^l}^*$ & $p^l-1$ & 0 & $\frac{m-1}{p^l}-d$\\
\hline
IIa & Affine $x_0\neq y_0$, $x_0$, $y_0\in\F_{p^l}^*$ & $p^l+1$ & $\left(\frac{p^l+1}{2}\right)^2$ & $(d-1)(d-2)$\\
\hline
IIb & Affine $x_0\neq y_0$, $x_0$ or $y_0\not\in\F_{p^l}^*$ & $p^l$ & 0 & $N_1$ \footnote{$N_1=\left(\frac{m-1}{p^l}-1\right)\left(2\frac{m-1}{p^l}-(m_b+1)p^{i_b-l}-1\right)-(d-1)(d-2)$}\\
\hline
IIc & Affine $x_0\neq y_0$, $x_0$ and $y_0\not\in\F_{p^l}^*$ & $p^l$ & $p^l$\footnote{$I_t(u,v)=0$ if $y_0(x_0+1)^{p^l}(y_0^{p^l-1}-1)^{p^l+1}\neq x_0(y_0+1)^{p^l}(x_0^{p^l-1}-1)^{p^l+1}$}  & $N_2$ \footnote{$N_2=\left\{\begin{array}{l}\left(\frac{m-1}{p^l}-1\right)\left(2\frac{m-1}{p^l}-(m_b+1)p^{i_b-l}-1\right)-(d-1)(d-2)\\\textrm{or }((p^l-2)(p^l+1)+1)(\frac{m-1}{p^l}-1) \\\hspace{1cm}\textrm{if $y_0(x_0+1)^{p^l}(y_0^{p^l-1}-1)^{p^l+1}= x_0(y_0+1)^{p^l}(x_0^{p^l-1}-1)^{p^l+1}$}\end{array}\right.$}\\
\hline
IIIa & $(1:1:0)$ & $p^l-1$ & $\left(\frac{p^l-1}{2}\right)^2$ & 1 \\
\hline
IIIb & $(\omega:1:0)$, \quad $\omega^d=1$ et $\omega\neq 1$ & $p^l$ & $\frac{p^{2l}-1}{4}$ & $d-1$\\
\hline
IIIc & $(\omega:1:0)$,\quad $\omega^d\neq1$ & $p^l-1$ &  0 & $\frac{m-1}{p^l}$ \\
\hline
\end{tabular}
\end{minipage}
\vspace{0.5cm}
\hrule
\end{table}

\subsection{Singular points at infinity}
We denote by $\widehat{f}$ (respectively $\widehat{h}$) the homogenized form of $f$ (respectively $h$). We denote by $\widetilde{f}$ (respectively $\widetilde{h}$) the dehomogenized form of $\widehat{f}$ (respectively $\widehat{h}$) relative to $y$.
\\\\Let $F(x,y,z)=(x+z)^m-x^m-(y+z)^m+y^m=z\widehat{f}$. Then, $$\left\{\begin{array}{rcl} F_x&=&m(x+z)^{m-1}-mx^{m-1}\\F_y&=&-m(y+z)^{m-1}+my^{m-1}\\F_z&=&m(x+z)^{m-1}-m(y+z)^{m-1}\end{array}\right..$$

At infinity ($z=0$), $F_x(x,y,0)=F_y(x,y,0)=0$ and $$F_z(x,y,0)=m(x^{m-1}-y^{m-1}).$$ So $(x_0,y_0,0)$ is a singular point of F if and only if $x_0^{m-1}=y_0^{m-1}$.
If $y_0=0$ then $x_0=0$; so $y_0\neq0$ and we have to study the solutions of \begin{equation}\label{1}x_0^{m-1}=1.\end{equation}
Equation \eqref{1} is equivalent to $x_0^{\frac{m-1}{p^l}}=1$. Since $\pgcd\left(\frac{m-1}{p^l},p\right)=1$, there are $\frac{m-1}{p^l}$ solutions at \eqref{1} and $x_0=1$ is the only one such that $x_0=y_0$.
\\\\\indent Now, we want to find the multiplicity of these singularities : \begin{align*}\widetilde{F}(x+x_0,z)&=(x+x_0+z)^m-(x+x_0)^m-(z+1)^m+1 \\&=\sum_{k=2}^m\binom{m}{k}(x+z)^kx_0^{m-k}-\sum_{k=2}^m\binom{m}{k}x^kx_0^{m-k}-\sum_{k=2}^m\binom{m}{k}z^k.\end{align*}
Since $m-1\equiv 0\mod p^l$, for all $2\leq k<p^l$, $\binom{m}{k}=0$.
Consider the terms of degree $p^l-1$ of $\widetilde{f}$ : $$\frac{1}{z}\binom{m}{p^l}(x_0^{m-p^l}(x+z)^{p^l}-x_0^{m-p^l}x^{p^l}-z^{p^l})=\binom{m}{p^l}(x_0^{m-p^l}-1)z^{p^l-1}.$$
This term vanishes (which means that $(x_0,y_0,0)$ is a singular point of multiplicity greater than $p^l-1$) if and only if $$x_0^{m-p^l}=1$$ that is to say if and only if $$ x_0^d=1.$$
Now, consider the terms of degree $p^l$ of $\widetilde{f}$ : \begin{align*}&\frac{1}{z}\binom{m}{p^l+1}(x_0^{m-p^l-1}(x+z)^{p^l+1}-x_0^{m-p^l-1}x^{p^l+1}-z^{p^l+1}) \\&\hspace{2.5cm}=\binom{m}{p^l+1}(x_0^{m-p^l-1}x^{p^l}+x_0^{m-p^l-1}xz^{p^l-1}+(x_0^{m-p^l-1}-1)z^{p^l}).\end{align*}
Since $x_0^{m-p^l-1}\neq0$, singular points of $\widehat{f}$ of multiplicity greater than  $p^l-1$ have multiplicity $p^l$.
\\\\We have just proved the following lemma :
\begin{lemme}\label{infini}
Let $\omega$ such that $\omega^{\frac{m-1}{p^l}}=1$. The point $(\omega:1:0)$ is a singular point of $\widehat{h}$ with multiplicity
$$\left\{\begin{array}{ll}p^l&\textrm{if $\omega^d=1$, $\omega\neq1$}\\p^l-1&\textrm{otherwise}\end{array}\right..$$
Furthermore, $\widehat{h}$ has $\frac{m-1}{p^l}$ singular points at infinity.\end{lemme}

\subsection{Affine singular points}
 We have : $$\left\{\begin{array}{rcl}f_x&=&m(x+1)^{m-1}-mx^{m-1}\\f_y&=&-m(y+1)^{m-1}+my^{m-1}\end{array}\right..$$
 So, \begin{align*}\textrm{$(x_0,y_0)$ singular point of $f$} &\Leftrightarrow\left\{\begin{array}{l}f(x_0,y_0)=0\\(x_0+1)^{m-1}=x_0^{m-1}\\(y_0+1)^{m-1}=y_0^{m-1}\end{array}\right. \\&\Leftrightarrow\left\{\begin{array}{l}x_0^{m-1}(x_0+1)-x_0^m-y_0^{m-1}(y_0+1)+y_0^m=0\\(x_0+1)^{m-1}=x_0^{m-1}\\(y_0+1)^{m-1}=y_0^{m-1}\end{array}\right. \\&\Leftrightarrow\left\{\begin{array}{l}x_0^{m-1}=y_0^{m-1}\\(x_0+1)^{m-1}=x_0^{m-1}\\(y_0+1)^{m-1}=y_0^{m-1}\end{array}\right..\end{align*}

 \begin{lemme}\label{affine}Affine singular points of $f$ are points satisfying $$(x_0+1)^{m-1}=x_0^{m-1}=y_0^{m-1}=(y_0+1)^{m-1}.$$\end{lemme}

 From Lemma \ref{affine}, we get that $x_0$, $y_0\neq0,-1$. Since $p^l$ divides $m-1$, \begin{equation}\label{2}\textrm{$(x_0,y_0)$ singular point of $f$}\Leftrightarrow \left\{\begin{array}{l}x_0^{\frac{m-1}{p^l}}=y_0^{\frac{m-1}{p^l}} \\(x_0+1)^{\frac{m-1}{p^l}}=x_0^{\frac{m-1}{p^l}} \\(y_0+1)^{\frac{m-1}{p^l}}=y_0^{\frac{m-1}{p^l}}\end{array}\right..\end{equation}
There are at most $\frac{m-1}{p^l}-1$ solutions to the second equation of \eqref{2}. Let $x_0$ be one of these solutions, we want to know the number of $y_0$ such that $(x_0,y_0)$ is a singular point of $f$.
\\\\\indent We write $m=1+\displaystyle\sum_{j=1}^b m_j p^{i_j}$ with $1\leq m_j\leq p-1$, $i_j>i_{j-1}$, $i_1=l$. Then,
\begin{align*}(y_0+1)^{\frac{m-1}{p^l}}=y_0^{\frac{m-1}{p^l}}&\Leftrightarrow\prod_{j=1}^b(y_0+1)^{m_jp^{i_j-l}}=y_0^{\frac{m-1}{p^l}} \\&\Leftrightarrow\sum_{0\leq k_j\leq m_j}^*\left(\prod_{j=1}^b\binom{m_j}{k_j}\right)y_0^{\sum_{j=1}^bk_jp^{i_j-l}}=0\end{align*}
where * indicates that this sum runs over all possible b-uples except $(m_1,\ldots,m_b)$.
We multiply by $y_0^{\frac{m-1}{p^l}-m_bp^{i_b-l}}$ and we set $\alpha=y_0^{\frac{m-1}{p^l}}$ : \begin{align*}&\sum_{\begin{subarray}{c}0\leq k_j\leq m_j\\j\neq b\end{subarray}}^*\left(\prod_{j=1}^{b-1}\binom{m_j}{k_j}\right)\alpha y_0^{\sum_{j=1}^{b-1}k_jp^{i_j-l}}\\&\hspace{2cm}+\sum_{k_b=0}^{m_b-1}\sum_{\begin{subarray}{c}0\leq k_j\leq m_j\\j\neq b\end{subarray}}\left(\prod_{j=1}^{b}\binom{m_j}{k_j}\right)y_0^{\frac{m-1}{p^l}-(m_b-k_b)p^{i_b-l}+\sum_{j=1}^{b-1}k_jp^{i_j-l}}=0.\end{align*}
The degree of this polynomial in $y_0$ is $$\frac{m-1}{p^l}-p^{i_b-l}+\displaystyle\sum_{j=1}^{b-1}m_jp^{i_j-l}=2\frac{m-1}{p^l}-(m_ b+1)p^{i_b-l}.$$
\begin{lemme}\label{nbreaff}The number of affine singularities of $h$ is at most :
$$\left(\frac{m-1}{p^l}-1\right)\left(2\frac{m-1}{p^l}-(m_b+1)p^{i_b-l}\right)$$ where $m=1+\displaystyle\sum_{j=1}^b m_j p^{i_j}$ with $1\leq m_j\leq p-1$, $i_j>i_{j-1}$, $i_1=l$.\end{lemme}

Now, we study the multiplicity of affine singularities :
\begin{align*}f(x+x_0,y+y_0)&=(x+x_0+1)^m-(x+x_0)^m-(y+y_0+1)^m+(y+y_0)^m \\&=\sum_{k=2}^m\binom{m}{k}x^k(x_0+1)^{m-k}-\sum_{k=2}^mx^kx_0^{m-k} \\&\hspace{2cm}-\sum_{k=2}^m\binom{m}{k}y^k(y_0+1)^{m-k}+\sum_{k=2}^my^ky_0^{m-k}.\end{align*}

Since $m-1\equiv0\mod p^l$, for all $2\leq k<p^l$, $\binom{m}{k}=0$. So $(x_0,y_0)$ is a singularity of multiplicity at least $p^l$. Consider the terms of degree $p^l+1$ :$$\binom{m}{p^l+1}(((x_0+1)^{m-p^l-1}-x_0^{m-p^l-1})x^{p^l+1}-((y_0+1)^{m-p^l-1}-y_0^{m-p^l-1})y^{p^l+1}).$$
Since $(x_0,y_0)$ is a singular point, $(x_0+1)^{m-1}=x_0^{m-1}$ and $x_0\neq-1,0$. So, \begin{align*}(x_0+1)^{m-p^l-1}-x_0^{m-p^l-1}=0&\Leftrightarrow (x_0+1)^{p^l}((x_0+1)^{m-p^l-1}-x_0^{m-p^l-1})=0\\&\Leftrightarrow -x_0^{m-p^l-1}=0.\end{align*}
Hence, affine singularities have multiplicity at most $p^l+1$.
Then we look at terms of degree $p^l$ :
$$\binom{m}{p^l}(((x_0+1)^{m-p^l}-x_0^{m-p^l})x^{p^l}-((y_0+1)^{m-p^l}-y_0^{m-p^l})y^{p^l}).$$
However,
\begin{align*}(x_0+1)^{m-p^l}-x_0^{m-p^l}=0&\Leftrightarrow(x_0+1)^{p^l}((x_0+1)^{m-p^l}-x_0^{m-p^l})=0 \\&\Leftrightarrow(x_0+1)^{m-1}(x_0+1)-x_0^m-x_0^{m-p^l}=0 \\&\Leftrightarrow x_0^{m-p^l}(x_0^{p^l-1}-1)=0\\&\Leftrightarrow x_0\in\F_{p^l}^*.\end{align*}
We can do the same for $y_0$.
\begin{lemme}\label{mulaff}There are at most :
\begin{itemize}\item $d-1$ affine singularities of $h$ such that $x_0=y_0\in\F_{p^l}^*$. They have multiplicity $p^l$ ($p^l+1$ for $f$);
\item $\frac{m-1}{p^l}-d$ affine singularities of $h$ such that $x_0=y_0\not\in\F_{p^l}^*$. They have multiplicity $p^l-1$ ($p^l$ for $f$);
\item$(d-1)(d-2)$ affine singularities of $h$ such that $x_0\neq y_0$ et $x_0$, $y_0\in\F_{p^l}^*$. They have multiplicity $p^l+1$ (for $h$ and $f$);
\item $\left(\frac{m-1}{p^l}-1\right)\left(2\frac{m-1}{p^l}-(m_b+1)p^{i_b-l}-1\right)-(d-1)(d-2)$ affine singularities of $h$ such that $x_0\neq y_0$ and $x_0$ or $y_0\not\in\F_{p^l}^*$. They have multiplicity $p^l$ (for $h$ and $f$).\end{itemize}\end{lemme}

\section{Intersection number bounds}\label{P3}

We write $h=uv$; we want to bound the intersection number $I_t(u,v)$ for $t$ a singularity of $h$. In the following, we use freely this lemma proved in \cite{MR1359909} :

\begin{lemme}\label{lemme1} Let $J(x,y)=0$ be an affine curve over $\F_q$ and $t=(x_0,y_0)$ be a point of $J$ of multiplicity $m_t$. Then $J(x+x_0,y+y_0)=J_{m_t}+J_{m_t+1}+\ldots$ where $J_i$ is an homogeneous polynomial of degree $i$; the factors of $J_m$ (on an algebraic closure) are called the tangent lines of $J$ at $t$. We write $J(x,y)=u(x,y).v(x,y)$; if $J_{m_t}$ and $J_{m_t+1}$ are relatively prime then $I_t(u,v)=m_t(u).m_t(v)$. Furthermore, if $J$ has only one tangent line at $t$, $I_t(u,v)=0$.\end{lemme}

\subsection{Singularities at infinity}

Let $t=(\omega:1:0)$ a singular point of $h$ at infinity ($\omega^{\frac{m-1}{p^l}}=1$). We write $\widetilde{h}(x+\omega,z)=\widetilde{H}_{m_t}+\widetilde{H}_{m_t+1}+\ldots$ where $m_t$ is the multiplicity of $t$ and $\widetilde{H}_i$ is the homogeneous polynomial composed of the terms of degree $i$ of $\widetilde{h}(x+\omega,z)$.
Then, \begin{align*}\widetilde{f}(x+\omega,z)&=\widetilde{h}(x+\omega,z)(x+\omega-1) \\&=(R+\widetilde{H}_{m_t+1}+\widetilde{H}_{m_t})(x+\omega-1)\\&\hspace{1cm}\textrm{where $R$ is a polynomial of degree greater than $m_t+1$} \\&=x R+(\omega-1)R+x\widetilde{H}_{m_t}+(\omega-1)\widetilde{H}_{m_t+1}+(\omega-1)\widetilde{H}_{m_t}.\end{align*}
So, \begin{itemize}\item if $\omega\neq1$, we have $\widetilde{F}_{m_t}=(\omega-1)\widetilde{H}_{m_t}$ and $\widetilde{F}_{m_t+1}=x\widetilde{H}_{m_t}+(\omega-1)\widetilde{H}_{m_t+1}$;
\item if $\omega=1$, $\widetilde{F}_{m_t+1}=x\widetilde{H}_{m_t}$.\end{itemize}

\begin{lemme}\label{III}If $t=(\omega:1:0)$, $\omega^{\frac{m-1}{p^l}}=1$, is a singular point at infinity of $h$ with multiplicity $m_t$ then
\begin{itemize}\item if $\omega\neq1$, $\widetilde{F}_{m_t}=(\omega-1)\widetilde{H}_{m_t}$ and $\widetilde{F}_{m_t+1}=x\widetilde{H}_{m_t}+(\omega-1)\widetilde{H}_{m_t+1}$;
\item if $\omega=1$, $\widetilde{F}_{m_t+1}=x\widetilde{H}_{m_t}$.\end{itemize}\end{lemme}

\begin{corollaire}If $t=(1:1:0)$ then $$I_t(u,v)\leq\left(\frac{p^l-1}{2}\right)^2.$$\end{corollaire}

\begin{preuve}If $t=(1:1:0)$ then its multiplicity is $p^l-1$. By Lemma \ref{III}, $$\widetilde{H}_{m_t}=a(x^{p^l-1}+z^{p^l-1}),$$ and all its factors are different. So $I_t(u,v)=m_t(u)m_t(v)$. We get the result since $m_t(u)+m_t(v)=p^l-1$.\end{preuve}

\begin{corollaire}If $t=(\omega:1:0)$ such that $\omega^d=1$, $\omega\neq1$ then $$I_t(u,v)\leq\frac{p^{2l}-1}{4}.$$\end{corollaire}

\begin{preuve}The multiplicity of $t$ is $p^l$. By Lemma \ref{III},  $$(\omega-1)\widetilde{H}_{p^l}=\widetilde{F}_{p^l}=x^{p^l}\omega^{m-p^l-1}+xz^{p^l-1}\omega^{m-p^l-1}+(\omega^{m-p^l-1}-1)z^{p^l}.$$
So all factors of $\widetilde{H}_{p^l}$ are simple and $I_t(u,v)=m_t(u)m_t(v)$. We get the result since $m_t(u)+m_t(v)=p^l$.\end{preuve}

\begin{corollaire} If $t=(\omega:1:0)$ with $\omega^{\frac{m-1}{p^l}}$, $\omega^d\neq1$, then $$I_t(u,v)=0.$$\end{corollaire}

\begin{preuve}The multiplicity of $t$ is $p^l-1$.
By Lemma \ref{III},$$(\omega-1)\widetilde{H}_{p^l-1}=\widetilde{F}_{p^l-1}=\alpha z^{p^l-1}$$ and $$\widetilde{F}_{p^l}=x\widetilde{H}_{p^l-1}+(\omega-1)\widetilde{H}_{p^l}=x^{p^l}\omega^{m-p^l-1}+xz^{p^l-1}+z^{p^l}(\omega^{m-p^l-1}-1).$$
So, $\pgcd(\widetilde{H}_{p^l},\widetilde{H}_{p^l-1})=\pgcd(\widetilde{F}_{p^l},\widetilde{F}_{p^l-1})=1$ and by Lemma \ref{lemme1}, $I_t(u,v)=0$.\end{preuve}

\subsection{Affine singularities}

We write $h(x+x_0,y+y_0)=H_{m_t}+H_{m_t+1}+\ldots$ where $m_t$ is the multiplicity of $t$ and $H_i$ is the homogeneous polynomial composed of the terms of degree $i$ of $h(x+x_0,y+y_0)$.
\\\\Let $t=(x_0,y_0)$ be an affine singular point of $h$ with multiplicity $m_t$ such that $x_0=y_0$. Then
\begin{align*}f(x+x_0,y+y_0)&=h(x+x_0,y+y_0)(x+x_0-y-y_0)\\&=(R+H_{m_t+1}+H_{m_t})(x-y) \\&\hspace{0.5cm}\textrm{where $R$ is a polynomial of degree greater than $m_t+1$}\\&=(x-y)R+(x-y)H_{m_t+1}+(x-y)H_{m_t}\\&=F_{m_t+1}+F_{m_t+2}+\ldots\end{align*}
So, $F_{m_t+2}=(x-y)H_{m_t+1}$ and  $F_{m_t+1}=(x-y)H_{m_t}$. Furthermore, for some $a$, $F_{m_t+1}=a(x^{m_t+1}-y^{m_t+1})$ (see proof of Lemma \ref{mulaff}).

\begin{lemme}\label{I=}If $t=(x_0,y_0)$ is an affine singular point of $h$ with multiplicity $m_t$ such that $x_0=y_0$, then
$F_{m_t+2}=(x-y)H_{m_t+1}$ and $F_{m_t+1}=(x-y)H_{m_t}$.
\\Furthermore, tangent lines to $h$ at $t$ are factors of $\frac{x^{m_t+1}-y^{m_t+1}}{x-y}$.\end{lemme}

\begin{corollaire}
Affine singularities of $h$, $t=(x_0,y_0)$, such that $x_0=y_0\in\F_{p^l}^*$ satisfy
$$I_t(u,v)\leq\frac{p^{2l}-1}{4}.$$\end{corollaire}

\begin{preuve}The multiplicity of $t$ is $p^l$. The factors of $\frac{x^{p^l+1}-y^{p^l+1}}{x-y}$ are all distinct. So, by Lemma \ref{I=}, tangent lines to $u$ or $v$ are all distinct and $I_t(u,v)=m_t(u)m_t(v)$. Since $m_t(u)+m_t(v)=p^l$, we get the result.\end{preuve}

\begin{corollaire}Affine singularities of $h$, $t=(x_0,y_0)$, such that $x_0=y_0\not\in\F_{p^l}^*$ satisfy $$I_t(u,v)=0.$$\end{corollaire}

\begin{preuve}The multiplicity of $t$ is $p^l-1$. By Lemma \ref{I=}, $$\textrm{$H_{p^l-1}=a(x-y)^{p^l-1}$ and $H_{p^l}=\frac{b(x^{p^l+1}-y^{p^l+1})}{x-y}$}.$$ Hence, $\pgcd(H_{p^l-1},H_{p^l})=1$. By Lemma \ref{lemme1}, $I_t(u,v)=0$.\end{preuve}

Let $t=(x_0,y_0)$ be an affine singular point of $h$ with multiplicity $m_t$ such that $x_0\neq y_0$. Then
\begin{align*}f(x+x_0,y+y_0)&=h(x+x_0,y+y_0)(x+x_0-y-y_0)\\&=(R+H_{m_t+1}+H_{m_t})(x+x_0-y-y_0)\\&\hspace{0.5cm}\textrm{where $R$ is a polynomial of degree greater than $m_t+1$}\\&=(x_0-y_0)H_{m_t}+((x-y)H_{m_t}+(x_0-y_0)H_{m_t+1})\\&\hspace{2cm}+((x-y+x_0-y_0)R+(x-y)H_{m_t+1})\\&=F_{m_t}+F_{m_t+1}+R' \\&\hspace{0.5cm}\textrm{where $R'$ is a polynomial of degree greater than $m_t+1$.}\end{align*}
So, $F_{m_t}=(x_0-y_0)H_{m_t}$ and $F_{m_t+1}=(x_0-y_0)H_{m_t+1}+(x-y)H_{m_t}$.

\begin{lemme}\label{II}If $t=(x_0,y_0)$ is an affine singular point of $h$ with multiplicity $m_t$ such that $x_0\neq y_0$ then
$$\textrm{$\F_{m_t}=(x_0-y_0)H_{m_t}$ and $F_{m_t+1}=(x-y)H_{m_t}+(x_0-y_0)H_{m_t+1}$}.$$\end{lemme}

\begin{corollaire} If $t=(x_0,y_0)$ is an affine singular point of $h$ such that $x_0\neq y_0$, $x_0$, $y_0\in\F_{p^l}^*$ then
$$I_t(u,v)\leq\left(\frac{p^l+1}{2}\right)^2.$$\end{corollaire}

\begin{preuve}The multiplicity of $t$ is $p^l+1$. By Lemma \ref{II}, $$\textrm{$(x_0-y_0)H_{m_t}=F_{m_t}=c_1x^{p^l+1}-c_2y^{p^l+1}$ with $c_1$, $c_2\neq 0$.}$$ Hence, all factors of $H_{m_t}$ are simple and then $I_t(u,v)=m_t(u)m_t(v)$. Since $m_t(u)+m_t(v)=p^l+1$, we get the result.\end{preuve}

\begin{corollaire}If $t=(x_0,y_0)$ is an affine singular point of
$h$ such that $x_0\neq y_0$ and $x_0\in\F_{p^l}^*$ and $y_0\not\in\F_{p^l}^*$ or $x_0\not\in\F_{p^l}^*$ and $y_0\in\F_{p^l}^*$ then $$I_t(u,v)=0.$$\end{corollaire}

\begin{preuve}The multiplicity of $t$ is $p^l$. Then
$$\textrm{$F_{p^l}=\left\{\begin{array}{ll}c_1x^{p^l}&\textrm{if $y_0\in\F_{p^l}^*$, $c_1\neq0$}\\c_2y^{p^l}&\textrm{if $x_0\in\F_{p^l}^*$, $c_2\neq0$}\end{array}\right.$
and $F_{p^l+1}=c_1'x^{p^l+1}-c_2'y^{p^l+1}$, $c_1'$, $c_2'\neq 0$.}$$
So, by Lemma \ref{II}, $1=\pgcd(F_{p^l},F_{p^l+1})=\pgcd(H_{p^l},H_{p^l+1})$.
Hence, by Lemma \ref{lemme1}, $I_t(u,v)=0$.\end{preuve}

Let $t=(x_0,y_0)$ be an affine singular point of $h$ such that $x_0\neq y_0$ and $x_0$, $y_0\not\in\F_{p^l}$; $t$ has multiplicity $p^l$.
We have $F_{p^l}=c_1x^{p^l}-c_2y^{p^l}=(c_3x-c_4y)^{p^l}$, where $c_1=(x_0+1)^{m-p^l}-x_0^{p^l}$, $c_2=(y_0+1)^{m-p^l}-y_0^{m-p^l}$; since $x_0$, $y_0\not\in \F_{p^l}^*$, $c_1\neq0$ and $c_2\neq 0$.
By Lemma \ref{II}, $$\textrm{$F_{p^l}=(x_0-y_0)H_{p^l}$ and $F_{p^l+1}=(x_0-y_0)H_{p^l+1}+(x-y)H_{p^l}$;}$$ so, $H_{p^l}$ has only one factor and $\pgcd(F_{p^l},F_{p^l+1})=\pgcd(H_{p^l},H_{p^l+1})$.
Furthermore, $F_{p^l+1}=d_1x^{p^l+1}-d_2y^{p^l+1}$ with $d_1=(x_0+1)^{m-p^l-1}-x_0^{m-p^l-1}\neq0$ and $d_2=(y_0+1)^{m-p^l-1}-y_0^{m-p^l-1}\neq0$.
Polynomials $F_{p^l}$ and  $F_{p^l+1}$ have a common factor if and only if $c_3x-c_4y$ divides $F_{p^l+1}$. So, $F_{p^l}$ and $F_{p^l+1}$ have a common factor if and only if $$\left(\frac{c_1}{c_2}\right)^{p^l+1}=\left(\frac{d_1}{d_2}\right)^{p^l}.$$
If $(x_0,y_0)$ is a singular point of $f$, then
$$\left\{\begin{array}{l}x_0^{m-1}=y_0^{m-1}\\(x_0+1)^{m-1}=x_0^{m-1}\\(y_0+1)^{m-1}=y_0^{m-1}\end{array}\right..$$
We have : \begin{align*}d_1=(x_0+1)^{m-p^l-1}-x_0^{m-p^l-1}&=\frac{(x_0+1)^{m-1}-x_0^{m-p^l-1}(x_0+1)^{p^l}}{(x_0+1)^{p^l}} \\&=\frac{x_0^{m-1}-x_0^{m-1}-x_0^{m-p^l-1}}{(x_0+1)^{p^l}}\\&=\frac{-x_0^{m-p^l-1}}{(x_0+1)^{p^l}}.\end{align*}
Similarly, $d_2=\frac{-y_0^{m-p^l-1}}{(y_0+1)^{p^l}}$.
Hence, $$\frac{d_1}{d_2}=\frac{x_0^{m-p^l-1}(y_0+1)^{p^l}}{y_0^{m-p^l-1}(x_0+1)^{p^l}} =\frac{x_0^{m-1}y_0^{p^l}(y_0+1)^{p^l}}{y_0^{m-1}x_0^{p^l}(x_0+1)^{p^l}}=\frac{y_0^{p^l}(y_0+1)^{p^l}}{x_0^{p^l}(x_0+1)^{p^l}}.$$
\\\\On the other hand, we have : \begin{align*}c_1=(x_0+1)^{m-p^l}-x_0^{m-p^l}&=\frac{(x_0+1)(x_0+1)^{m-1}-x_0^{m-p^l}(x_0+1)^{p^l}}{(x_0+1)^{p^l}} \\&=\frac{x_0^{m}+x_0^{m-1}-x_0^{m}-x_0^{m-p^l}}{(x_0+1)^{p^l}}\\&=\frac{x_0^{m-p^l}(x_0^{p^l-1}-1)}{(x_0+1)^{p^l}}.\end{align*}
Similarly, $c_2=\frac{y_0^{m-p^l}(y_0^{p^l-1}-1)}{(y_0+1)^{p^l}}$.
Hence, $$\frac{c_1}{c_2}=\frac{x_0^{m-p^l}(x_0^{p^l-1}-1)(y_0+1)^{p^l}}{y_0^{m-p^l}(y_0^{p^l-1}-1)(x_0+1)^{p^l}} =\frac{y_0^{p^l-1}(y_0+1)^{p^l}(x_0^{p^l-1}-1)}{x_0^{p^l-1}(x_0+1)^{p^l}(y_0^{p^l-1}-1)}.$$
\indent After simplification, we get that $F_{p^l}$ and $F_{p^l+1}$ have a common factor if and only if \begin{equation}\label{*}y_0(x_0+1)^{p^l}(y_0^{p^l-1}-1)^{p^l+1}=x_0(y_0+1)^{p^l}(x_0^{p^l-1}-1)^{p^l+1}.\end{equation}
\indent If $(x_0,y_0)$ is not a solution of \eqref{*}, then $\pgcd(H_{p^l},H_{p^l+1})=1$ and by Lemma \ref{lemme1}, $I_t(u,v)=0$.
\\\\\indent Otherwise, we write $u(x+x_0,y+y_0)=U_r+U_{r+1}+\ldots$, with $U_r\neq0$
and $v(x+x_0,y+y_0)=V_s+V_{s+1}+\ldots$, with $V_s\neq0$.
If $r=0$ or $s=0$ then $t$ is not a point of $u$ or $v$ and $I_t(u,v)=0$.
Assume that $r$, $s>0$. Since $(x_0,y_0)$ satisfy \eqref{*}, $F_{p^l}$ and $F_{p^l+1}$ have a common factor that we denote by $e$. We have $H_{p^l}=U_rV_s=e^{p^l}$ and $H_{p^l+1}=U_rV_{s+1}+U_{r+1}V_s$. Furthermore, $\pgcd(F_{p^l},F_{p^l+1})=e$ and thus $\pgcd(H_{p^l},H_{p^l+1})=e$. Since $r\geq1$ and $s\geq1$, $e$ divides $U_r$ and $V_s$ and consequently $\pgcd(U_r,V_s)$. If $\pgcd(U_r,V_s)=e^k$, $e^k$ divides $\pgcd(H_{p^l},H_{p^l+1})$ thus $\pgcd(U_r,V_s)=e$.
We can assume without loss of generality that $U_r=e^{p^l-1}$ and $V_s=e$; $t$ is a simple point of $v$ thus $I_t(u,v)=\mathrm{ord}_t^v(u)$ (see \cite{fulton_alg_curves}). Since $e^2$ does not divide $H_{p^l+1}$, $e$ does not divide $U_{p^l}$ and we can write $U_{p^l}$ as the product of $p^l$ linear factors distinct from $e$. Each factor is not tangent to $v$, so the order of each factor is 1. Thus the order of $U_{p^l}$ is $p^l$ and $\mathrm{ord}_t^v(u)\leq p^l$.

\begin{lemme}\label{IIc}If $t=(x_0,y_0)$ is an affine singular point of $h$ such that $x_0\neq y_0$ and $x_0$ and $y_0\not\in \F_{p^l}^*$
then \begin{itemize}\item if $y_0(x_0+1)^{p^l}(y_0^{p^l-1}-1)^{p^l+1}\neq x_0(y_0+1)^{p^l}(x_0^{p^l-1}-1)^{p^l+1}$, $I_t(u,v)=0$
\item otherwise, $I_t(u,v)\leq p^l$;
and there is at most $((p^l-2)(p^l+1)-1)(\frac{m-1}{p^l}-1)$ such singular points.\end{itemize}\end{lemme}

\section{Proof of theorem \ref{facteurs}}

\begin{lemme}\label{cle}If $h$ has no absolutely irreducible factor over $\F_p$ then there exists a factorization $h=uv$ such that $$\sum_{t}I_t(u,v)\geq 2\frac{\mathrm{deg}(h)^2}{9}.$$
Equivalently, if $I_{tot}$ is any upper bound on the global intersection number $\displaystyle\sum_tI_t(u,v)$ of $u$ and $v$ for all factorizations of $h$ into two factors over the algebraic closure of $\F_p$, then $$e=\frac{I_{tot}}{\frac{\mathrm{deg}(h)^2}{4}}\geq\frac{8}{9}.$$\end{lemme}

\begin{preuve}We write $h=e_1\ldots e_r$, where each $e_i$ is irreducible over $\F_p$, but not absolutely irreducible. Then each $e_i$ factors into $c_i\geq2$ factors on an algebraic closure of $\F_p$ and its factors are all of degree $\frac{\mathrm{deg}(e_i)}{c_i}$.
Now, we factor each $e_i$ into two factors $u_i$ and $v_i$ such that $\mathrm{deg}(u_i)=\mathrm{deg}(v_i)+\frac{\mathrm{deg}(e_i)}{c_i}$ if $c_i$ is odd (thus $c_i\geq3$) or  $\mathrm{deg}(u_i)=\mathrm{deg}(v_i)$ if $c_i$ is even.
We set $u=\displaystyle\prod_{i=1}^ru_i$ and $v=\displaystyle\prod_{i=1}^rv_i$.
Then $\deg(u)-\deg(v)\leq\frac{\deg(h)}{3}$. Since $\deg(u)+\deg(v)=\deg(h)$, $$\deg(u)\deg(v)\geq \frac{8}{9}\frac{\deg(h)^2}{4}.$$
Let $I_{tot}$ be an upper bound on the global intersection number of $u$ and $v$ for all factorizations of $h$ into two factors over the algebraic closure of $\F_p$. Then by Bezout's theorem, $$I_{tot}\geq\sum_tI_t(u,v)= \deg(u)\deg(v)\geq\frac{8}{9}\frac{\deg{(h)}^2}{4}=2\frac{\mathrm{deg}(h)^2}{9}.$$
\end{preuve}

The following lemma is proved in \cite{mcguire_excep} for $p=2$ but it is the exact same proof for $p\neq 2$.
\begin{lemme}\label{17}Let $h_k$, $1\leq k\leq r$, the irreducible factors of $h$ over $\F_p$ and for all $1\leq k\leq r$, we write $h_k=h_{k,1}\ldots h_{k,c_k}$ the factorization of $h_k$ into $c_k$ absolutely irreducible factors. Then \begin{enumerate}\item $\deg(h_k)^2\leq \displaystyle\sum_{t \in Sing(F)}m_t(h_k)^2$ where $Sing(F)$ is the set of singular points of $F$.
\item $\displaystyle\sum_{1\leq i< j\leq c_k}m_t(h_{k,i})m_t(h_{k,j})\leq m_t(h_k)^2\frac{c_k-1}{2c_k}$. \end{enumerate}\end{lemme}

The following theorems prove Theorem \ref{facteurs}. From now, we assume $m\neq1+p^l$.

\begin{theoreme}If $d=1$ then $h$ has an absolutely irreducible factor over $\F_p$.\end{theoreme}

\begin{preuve}Assume that $h$ has no absolutely irreducible factor over $\F_p$, then by Lemma \ref{cle} we must have $e=\frac{I_{tot}}{\frac{(m-2)^2}{4}}\geq\frac{8}{9}$ where $I_{tot}$ is an upper bound on the global intersection number.
Since $d=1$, we only have singularities of type Ib, IIc, IIIa et IIIc (see Table \ref{tab1}). So
\begin{equation}\label{eq1}\sum_{t}I_t(u,v)\leq p^l\left(\frac{m-1}{p^l}-1\right)\left(2\frac{m-1}{p^l}-(m_b+1)p^{i_b-l}-1\right)+\left(\frac{p^l-1}{2}\right)^2.\end{equation}
Since $m=1+p^lk$ and $m\neq1+p^l$, $k\geq2$; thus  $\frac{m-3}{4}=\frac{p^lk-2}{4}\geq \frac{p^l-1}{2}$. Hence
\begin{align*}e&\leq\frac{1}{\frac{(m-2)^2}{4}}\left(\frac{(m-3)^2}{16}+p^l(\frac{m-1}{p^l}-1)^2\right)\\&\leq\frac{1}{4}+\frac{4}{p^l}.\end{align*}
For $p^l\neq3$ or 5, we have $e<\frac{8}{9}$ which is a contradiction.
\\\\First, consider the case where $p^l=3$. We have $1=d=\pgcd(2,k)$ so $k$ is odd and 3 does not divide $k$ by definition of $l$. Hence $k\geq5$, thus \begin{align*}e&\leq\frac{p^l((p^l-2)(p^l+1)+1)(\frac{m-1}{p^l}-1)+\left(\frac{p^l-1}{2}\right)^2}{\frac{(m-2)^2}{4}}=\frac{15(k-1)+1}{\frac{(3k-1)^2}{4}}.\end{align*} However, for $k\geq 5$, $k\mapsto \frac{15(k-1)+1}{\frac{(3k-1)^2}{4}}$ is a decreasing function. So, for $k\geq 11$, $e<8/9$. Now we have to consider the case where $k=5$ and $k=7$.  Using equation \eqref{eq1}, we have $$\begin{array}{|c|c|c|}\hline
k&5&7\\ \hline m&16&22\\ \hline I_{tot}&37&73 \\ \hline e&\frac{37}{7^2}&\frac{73}{11^2}\\ \hline\end{array}$$In all cases we get a contradiction.
\\\\\indent If $p^l=5$, $1=d=\pgcd(4,k)$ and $k$ is odd. Hence, $k=3$ or $k\geq7$.
As in the case where $p^l=3$, $e\leq\frac{95(k-1)+4}{\frac{(5k-1)^2}{4}}$.
However $k\mapsto\frac{95(k-1)+4}{\frac{(5k-1)^2}{4}}$ is a decreasing function for $k\geq3$. so, for $k\geq17$, $e<\frac{8}{9}$ which is a contradiction.
We now have to consider the case where $k=3$, 7, 9, 11, 13.
Using equation \eqref{eq1}, we have $$\begin{array}{|c|c|c|c|c|c|}\hline
k&3&7&9&11&13\\ \hline m&16&36&46&56&66\\ \hline I_{tot}&24&124&324&354&664 \\ \hline e&\frac{24}{7^2}&\frac{124}{17^2}&\frac{324}{22^2}&\frac{354}{27^2}&\frac{664}{32^2}\\ \hline\end{array}$$
In all case, $e<\frac{8}{9}$ which is a contradiction since $e<\frac{8}{9}$.\end{preuve}

\begin{theoreme}If $1<d<\frac{m-1}{p^l}$, $h$ has an absolutely irreducible factor over $\F_p$.\end{theoreme}

\begin{preuve}Assume that $h$ has no absolutely irreducible factor over $\F_p$, then by Lemma \ref{cle}, we must have $e=\frac{I_{tot}}{\frac{(m-2)^2}{4}}\geq\frac{8}{9}$ where $I_{tot}$ is an upper bound on the global intersection number.
We have : \begin{align*}\sum_tI(u,v)&\leq\frac{p^{2l}-1}{4}(d-1)+\left(\frac{p^l-1}{2}\right)^2  \\&+p^l\left((\frac{m-1}{p^l}-1)(2\frac{m-1}{p^l}-(m_b+1)p^{i_b-l}-1)-(d-1)(d-2)\right) \\&\hspace{3cm}+\left(\frac{p^l+1}{2}\right)^2(d-1)(d-2)+(d-1)\frac{p^{2l}-1}{4} \\&\leq\frac{p^{2l}-1}{2}(d-1)+\left(\frac{p^l-1}{2}\right)^2(d-1)(d-2) \\&\hspace{3cm}+p^l\left(\frac{m-1}{p^l}-1\right)^2+\left(\frac{p^l-1}{2}\right)^2.\end{align*}
However, $m=1+kp^l$ with $k\neq1$. Since $d$ divide $k$ and $d<k$ we have $d\leq \frac{m-1}{2p^l}$. Hence,
\begin{align*}e&\leq\frac{2(p^{2l}-1)(\frac{k}{2}-1)+(p^l-1)^2(\frac{k}{2}-1)(\frac{k}{2}-2)+4p^l(k-1)^2+(p^l-1)^2}{(p^lk-1)^2} \\&\leq\frac{1}{\left(k-\frac{1}{p^l}\right)^2}\left((1-\frac{1}{p^{2l}})(k-2)+\frac{1}{4}(1-\frac{1}{p^l})^2(k-2)(k-4)\right. \\&\hspace{6.5cm}+\left.\frac{4}{p^l}(k-1)^2+(1-\frac{1}{p^l})^2\right) \\e&\leq\frac{1}{k-\frac{1}{p^l}}+\frac{1}{4}+\frac{4}{p^l}+\frac{1}{\left(k-\frac{1}{p^l}\right)^2}.\end{align*}
Since $e\geq\frac{8}{9}$, $1<d<k$ and $\pgcd(k,p)=1$, the only possibilities are :
$$\begin{array}{|c|c|c|c|c|c|c|c|c|c|}\hline k&4&6&8&9&10&12&14&15&\geq16\\ \hline p^l&3,7,11&5&3,5,7&7&3,7&5,7&3,5&7&3,5\\ \hline\end{array}$$
On one hand, we have \begin{align}\label{eqa}e&\leq \frac{2(p^{2l}-1)(d-1)+(p^l+1)^2(d-1)(d-2)}{(p^lk-1)^2}\\ \nonumber&\hspace{2cm}+\frac{4p^l(k-1)\left((p^l-2)(p^l+1)+1\right)+(p^l-1)^2}{(p^lk-1)^2}.\end{align}
On the other hand, we have :
\begin{align}\label{eqb}e&\leq \frac{2(p^{2l}-1)(d-1)+(p^l+1)^2(d-1)(d-2)}{(p^lk-1)^2}\\ \nonumber&\hspace{2cm}+\frac{4p^l(k-1)(2k-(m_b+1)p^{i_b-l}-1)+(p^l-1)^2}{(p^lk-1)^2}.\end{align}
First, consider the case where $k\geq16$. In inequality \eqref{eqa}, e is bounded by a decreasing function of $k$. Furthermore, if $p^l=3$ and $k=16$ or if $k=17$ and $p^l=5$ the upper bound in \eqref{eqa} is less than $\frac{8}{9}$ which leaves only the case $k=16$ and $p^l=5$. But replacing in equation \eqref{eqb}, we also get a contradiction.
In the other cases, using inequality \eqref{eqa} or inequality \eqref{eqb}, we have $e<\frac{8}{9}$ which is a contradiction.\end{preuve}

\begin{theoreme}\label{t3} If $d=\frac{m-1}{p^l}\neq p^l-1$ then $h$ has an absolutely irreducible factor over $\F_p$.\end{theoreme}

\begin{preuve}First, we make some remarks; since $d=\frac{m-1}{p^l}$, there are only singularities of type Ia, IIa, IIIa, IIIb (see Table \ref{tab1}).
In all these case, $H_{m_t}$ only has simple factors. So, for all factorization $h=uv$, $I_t(u,v)=m_t(u)m_t(v)$. Furthermore, since $\frac{m-1}{p^l}\neq p^l-1$, $\frac{m-1}{p^l}\leq \frac{p^l-1}{2}$.
Assume that $h$ has no absolutely irreducible factor over $\F_p$. We write $h=h_1\ldots h_r$ where each $h_i$ factorizes into $c_i\geq2$ factors on an algebraic closure of $\F_p$ and its factors are all of degree $\frac{\mathrm{deg}(h_i)}{c_i}$. We write $h_i=h_{i,1}\ldots h_{i,c_i}$. Then \begin{align*}A&=\sum_{k=1}^r\sum_{1\leq i<j\leq c_k}\sum_t I_t(h_{k,i},h_{k,j})+\sum_{1\leq k<l\leq r}\sum_{\begin{subarray}{c}1\leq i\leq c_k\\1\leq j\leq c_l\end{subarray}}\sum_t I_t(h_{k,i},h_{l,j})\\&=\sum_{k=1}^r\sum_{1\leq i<j\leq c_k}\sum_t m_t(h_{k,i})m_t(h_{k,j})+\sum_{1\leq k<l\leq r}\sum_{\begin{subarray}{c}1\leq i\leq c_k\\1\leq j\leq c_l\end{subarray}}\sum_t m_t(h_{k,i})m_t(h_{l,j}).\end{align*}
However, \begin{align*}(m_t(h))^2&=\left(\sum_{k=1}^rm_t(h_k)\right)^2\\&=\sum_{k=1}^rm_t(h_k)^2+2\sum_{1\leq k<l\leq r}m_t(h_k)m_t(h_l) \\&=\sum_{k=1}^rm_t(h_k)^2+2\sum_{1\leq k<l\leq r}\sum_{\begin{subarray}{c}1\leq i\leq c_k\\1\leq j\leq c_l\end{subarray}} m_t(h_{k,i})m_t(h_{l,j}).\end{align*}
So, by Lemma \ref{17}, 
$$A\leq \sum_t\left( \sum_{k=1}^rm_t(h_k)^2\frac{c_k-1}{2c_k}+\frac{1}{2}(m_t(h)^2-\sum_{k=1}^rm_t(h_k)^2)\right),$$
thus
$$A\leq\frac{1}{2}\sum_t\left(m_t(h)^2-\sum_{k=1}^r\frac{m_t(h_k)^2}{c_k}\right).$$
On the other hand, by Bezout's Theorem, 
\begin{align*}A&=\sum_{k=1}^r\sum_{1\leq i<j\leq c_k}\deg(h_{k,i})\deg(h_{k,j})+\sum_{1\leq k<l\leq r}\sum_{\begin{subarray}{c}1\leq i\leq c_k\\1\leq j\leq c_l\end{subarray}}\deg(h_{k,i})\deg(h_{l,j})\\&=\sum_{k=1}^r\frac{\deg(h_k)^2}{c_k^2}\frac{c_k(c_k-1)}{2}+\sum_{1\leq k<l\leq r}\deg(h_k)\deg(h_l) \\&=\sum_{k=1}^r\deg(h_k)^2\frac{c_k-1}{2c_k}+\frac{1}{2}\left(\deg(h)^2-\sum_{k=1}^r\deg(h_k)^2\right) \\&=\frac{1}{2}\left(\deg(h)^2-\sum_{k=1}^r\frac{\deg(h_k)^2}{c_k}\right).\end{align*}
Hence, $$\deg(h)^2-\sum_{k=1}^r\frac{\deg(h_k)^2}{c_k}\leq \sum_t\left(m_t(h)^2-\sum_{k=1}^r\frac{m_t(h_k)^2}{c_k}\right).$$
Then, by Lemma \ref{17}, $$\deg(h)^2-\sum_tm_t(h)^2\leq \sum_{k=1}^r\frac{1}{c_k}\left(\deg(h_k)^2-\sum_tm_t(h_k)^2\right)\leq 0.$$
We set $k=\frac{m-1}{p^l}$. Then
\begin{align*}\deg(h)^2\leq\sum_tm_t(h)^2&\Leftrightarrow (m-2)^2\leq 2(k-1)p^{2l}\\&\hspace{2.5cm}+(k-1)(k-2)(1+p^l)^2+(p^l-1)^2 \\&\Leftrightarrow-(2p^l+1)k^2+(p^{2l}+4p^l+3)k-(p^{2l}+2p^l+2)\leq 0 \\&\Leftrightarrow \textrm{$k\leq 1$ ou $k\geq\frac{p^{2l}+2p^l+2}{2p^l+1}$}.\end{align*}
However, $k\geq2$ ($m\neq1+p^l$) and $k\leq\frac{p^l-1}{2}<\frac{p^{2l}+2p^l+2}{2p^l+1}$ which is a contradiction.\end{preuve}

\section{The case where $m\not\equiv1\mod p$}
In this section we assume that $m\not\equiv 0\mod p$ and $m\not\equiv 1\mod p$.
Studying singularities as in Parts \ref{P2} and \ref{P3}, we get the following lemma :

\begin{lemme}Affine singularities of $h$ are the points $(x_0,y_0)$ such that $$(x_0+1)^{m-1}=(y_0+1)^{m-1}=(x_0)^{m-1}=(y_0)^{m-1}$$ and $x_0\neq y_0.$ They have multiplicity 2. There is no singularity at infinity.\end{lemme}

Using the same kind of idea that in the proof of Theorem \ref{t3}, we get that if $h$ has no absolutely irreducible factor then $$(m-2)^2\leq 4N$$ where $N$ is the number of affine singularities.
\\\\\indent Unfortunately, we cannot eliminate any $m$ with this method. However, we have the following result :

\begin{proposition}If $\gcd(m-1,p-1)>1$ then $h$ has an absolutely irreducible factor.\end{proposition}

\begin{preuve}We write $h=ag_1\ldots g_r$ where $g_i$ are monic in $x$ and absolutely irreducible on an extension $E$ of $\F_p$. Let $h_i$ be the homogeneous polynomial of highest degree in $g_i$. Then $$\frac{x^{m-1}-y^{m-1}}{x-y}=h_1\ldots h_r.$$
On the other hand, since $m\not\equiv1\mod p$, $$\frac{x^{m-1}-y^{m-1}}{x-y}=\prod_{\zeta}(x-\zeta y)$$ where the product runs over all $(m-1)$th roots of unity but 1. Each $(x-\zeta y)$ divides only one of the $h_i$.
Now, we consider : $$\sigma : \begin{array}{ccc}E[x,y]&\rightarrow&E[x,y]\\\displaystyle\sum_{i,j}a_{i,j}x^iy^j&\mapsto&\displaystyle\sum_{i,j}a_{i,j}^px^iy^j\end{array}.$$
Then $h=\sigma(h)=a\sigma(g_1)\ldots\sigma(g_r)$. By unicity of factorization, for all $1\leq i\leq r$, there exists $1\leq j\leq r$ such that $\sigma(g_i)=g_j$ and consequently $\sigma(h_i)=h_j$. However, since $\gcd(m-1,p-1)>1$, there exists a $(m-1)$th root of unity $\zeta_0\neq1$ such that $\zeta_0\in\F_p$. So $\sigma(x-\zeta_0 y)=x-\zeta_0 y$. Assume for example that $x-\zeta_0 y$ divides $h_1$, then since $x-\zeta_0 y$ divides only one of the $h_i$, $\sigma(h_1)=h_1$. Finally, $\sigma(g_1)=g_1$ which means that $g_1\in\F_p[x,y]$ and $h$ as an absolutely irreducible  factor over $\F_p$. \end{preuve}

\bibliographystyle{plain}
\bibliography{C:/Users/Elodie/Dropbox/modele-these-bdd/bibliothese}

\end{document}